\documentclass[11pt]{article}
\usepackage[a4paper,top=3cm,bottom=2cm,left=2cm,right=2cm,marginparwidth=1.75cm]{geometry}

\usepackage[utf8]{inputenc} \usepackage[T1]{fontenc}
\usepackage[numbers]{natbib}
\usepackage[colorlinks=True,citecolor=blue,urlcolor=blue,pagebackref=true,backref=true]
{hyperref}
\renewcommand*{\backrefalt}[4]{\ifcase #1 \footnotesize{(Not cited.)}\or        \footnotesize{(Cited on page~#2.)}\else      \footnotesize{(Cited on pages~#2.)}\fi}

\usepackage[none]{hyphenat}
\usepackage{breakcites}
\usepackage[utf8]{inputenc} \usepackage[T1]{fontenc}    \usepackage{hyperref}       \hypersetup{
  colorlinks = true,
  urlcolor = blue,
  linkcolor = blue,
  citecolor = blue
}
\usepackage{booktabs}       \usepackage{amsfonts}       \usepackage{nicefrac}       \usepackage{microtype}      \usepackage{xcolor}         \usepackage{caption}
\usepackage{subcaption}
\usepackage{float}
\usepackage{amsmath,amsthm}
\usepackage{amssymb}
\usepackage{mathtools}
\usepackage{natbib}
\usepackage{soul}
\usepackage{bm}
\usepackage{enumitem}
\usepackage{multirow}
\usepackage{wrapfig}
\usepackage{scalerel}
\allowdisplaybreaks
\usepackage{dsfont}
\usepackage[noend]{algpseudocode}
\usepackage{xcolor}
\usepackage{thmtools,thm-restate}
\usepackage{cleveref}

 \newtheorem{theorem}{Theorem}
 \newtheorem{definition}{Definition}
 \newtheorem{lemma}{Lemma}
 
\newtheorem{proposition}{Proposition}

\title{Sub-optimality of the Separation Principle \\ for Quadratic Control from Bilinear Observations}
\author{Yahya Sattar \\ Cornell \\ {\tt ysattar@cornell.edu} \and \and Sunmook Choi \\ Cornell \\ { \tt sc3377@cornell.edu} \and \and  Yassir Jedra \\ MIT  \\ { \tt jedra@mit.edu} \and \and  Maryam Fazel \\ U Washington \\ { \tt mfazel@uw.edu} \and \and Sarah Dean \\ Cornell  \\ { \tt sdean@cornell.edu}}

\date{}
\usepackage{enumitem}
\usepackage{caption,subcaption,multirow}
\usepackage{bm}

\renewcommand{\cite}{\citep}

\newtheorem{assumption}{Assumption}

\newcommand{\bv}[1]{{\boldsymbol{#1}}}		\newcommand{\bvgrk}[1]{{\boldsymbol{#1}}}

\newcommand{\xhat}{\hat{x}}

\newcommand{\vu}{\bv{u}}
\newcommand{\vv}{\bv{v}}
\newcommand{\vw}{\bv{w}}
\newcommand{\vx}{\bv{x}}
\newcommand{\vy}{\bv{y}}
\newcommand{\vz}{\bv{z}}

\newcommand{\ub}{\bv{u}}
\newcommand{\wb}{\bv{w}}
\newcommand{\xb}{\bv{x}}

\newcommand{\vxhat}{\hat{\vx}}

\newcommand{\Acal}{\mathcal{A}}
\newcommand{\Bcal}{\mathcal{B}}

\newcommand{\Gcal}{\mathcal{G}}

\newcommand{\Ical}{\mathcal{I}}

\newcommand{\Ncal}{\mathcal{N}}

\newcommand{\vA}{\bv{A}}
\newcommand{\vB}{\bv{B}}
\newcommand{\vC}{\bv{C}}

\newcommand{\vI}{\bv{I}}

\newcommand{\vK}{\bv{K}}
\newcommand{\vL}{\bv{L}}

\newcommand{\vO}{\bv{O}}
\newcommand{\vP}{\bv{P}}
\newcommand{\vQ}{\bv{Q}}
\newcommand{\vR}{\bv{R}}

\newcommand{\vU}{\bv{U}}
\newcommand{\vV}{\bv{V}}

\newcommand{\vX}{\bv{X}}

\newcommand{\Ab}{\bv{A}}

\newcommand{\Bb}{\bv{B}}

\newcommand{\vSigma}{\bvgrk{\Sigma}}

\newcommand{\nn}{\nonumber}

\newcommand{\E}{\operatorname{\mathbb{E}}}

\newcommand{\eqsym}[1]{\stackrel{\text{(#1)}}{=}}

\newcommand{\norm}[1]{\|{#1} \|}

\newcommand{\R}{\mathbb{R}}				

\newcommand{\dm}[2]
{
	\IfStrEq{#2}{1}{\R^{#1}}{\R^{#1 \x #2}}
}

 			\newcommand{\tr}{\textup{\textbf{tr}}} 			 		 								 	 	  				 				 											\newcommand{\distas}{\overset{\text{i.i.d.}}{\sim}}

\newcommand{\T}{\top}
\newcommand*{\x}{\mathsf{x}\mskip1mu}

\newcommand{\splitatcommas}[1]{\begingroup
	\begingroup\lccode`~=`, \lowercase{\endgroup
		\edef~{\mathchar\the\mathcode`, \penalty0 \noexpand\hspace{0pt plus .1em}}}\mathcode`,="8000 #1\endgroup
}

\usepackage{bbm}
\newcommand{\Item}[1]{\ifx\relax#1\relax  \item \else \item[#1] \fi
	\abovedisplayskip=0pt\abovedisplayshortskip=0pt~\vspace*{-\baselineskip}
}

\usepackage{booktabs}
\usepackage{multirow}

\usepackage{mathtools}
\makeatletter
\newcommand{\raisemath}[1]{\mathpalette{\raisem@th{#1}}}\newcommand{\raisem@th}[3]{\raisebox{#1}{$#2#3$}}
\makeatother

\newcommand{\newcustomtheorem}[2]{\newenvironment{#1}[1]
	{\renewcommand\customgenericname{#2}\renewcommand\theinnercustomgeneric{##1}\innercustomgeneric
	}
	{\endinnercustomgeneric}
}
\newcustomtheorem{customtheorem}{Theorem}
\newcustomtheorem{customlemma}{Lemma}
\newcustomtheorem{customproposition}{Proposition}

\newcommand{\vertiii}[1]{{\vert\kern-0.25ex\vert\kern-0.25ex\vert #1 \vert\kern-0.25ex\vert\kern-0.25ex\vert}}

\newcommand{\beq}{\begin{equation}}
	
	\newcommand{\eeq}{\end{equation}}

\newcommand{\Ac}{{\cal{A}}}

\newcommand{\Gc}{{\cal{G}}}

\newcommand{\Bc}{\mathcal{B}}

\newcommand{\ubar}{\bar{u}}
\newcommand{\Lbar}{\bar{L}}
\newcommand{\dbar}{\bar{d}}

\definecolor{emmanuel}{RGB}{255,127,0}

\numberwithin{equation}{section}

\newcommand{\sd}[1]{\textcolor{teal}{[sd: #1]}}
\newcommand{\yj}[1]{\textcolor{purple}{[yj: #1]}}

\newcommand{\mf}[1]{\textcolor{magenta}{[mf: #1]}}

\newcommand{\fbar}{\bar{f}}

 \usepackage[]{mdframed}

\begin{document}

\maketitle

\begin{abstract}
We consider the problem of controlling a linear dynamical system from bilinear observations with minimal quadratic cost. 
Despite the similarity of this problem to standard linear quadratic Gaussian (LQG) control, we show that when the observation model is bilinear, neither does the Separation Principle hold, nor is the optimal controller affine in the estimated state.
Moreover, the cost-to-go is non-convex in the control input.
Hence, finding an analytical expression for the optimal feedback controller is difficult in general.  
Under certain settings, we show that the standard LQG controller locally maximizes the cost instead of minimizing it.
Furthermore, the optimal controllers~(derived analytically) are not unique and are nonlinear in the estimated state.
We also introduce a notion of input-dependent observability and derive conditions under which the Kalman filter covariance remains bounded. 
We illustrate our theoretical results through numerical experiments in multiple synthetic settings. 
\end{abstract}

\section{Introduction}
In many engineering settings, measurements result from the interaction between a measurement device and an unknown quantity.
Traditionally, the unknown quantity is assumed to be independent of the measurement process;
however, in general it may be influenced by being measured.
This ``observer effect'' is present in examples ranging from electronic circuits, where measurement devices can alter resistance or impedance~\cite{liu2025probabilistic}, to robotics, where active perception requires interaction~\cite{jamali2016active}, and quantum systems, where measurement induces wavefunction collapse~\cite{pardalos2010optimization}.
Perhaps the simplest model which captures such interaction is a \emph{bilinear} observation model,
in which the output is a bilinear function of the control input and the unknown state. Recent works~\cite{liu2025probabilistic,sattar2024learning} propose a dynamical system model with linear transitions and bilinear measurements, and study the system identification problem.

In such settings, control inputs both change the state and affect observation of it.
How should controllers be designed to account for these dual roles?
In this paper, we address this question for optimal control of a linear dynamical system with bilinear observations and quadratic costs.
The bilinearity is a small departure from the standard partially observed linear quadratic problem \cite{bertsekas2012dynamic} (LQG),
in which the optimal controller is known to obey the \emph{Separation Principle}: it suffices to independently consider state estimation and a state feedback controller. 
Separation principle has also been shown to  holds in several variants of LQG, e.g., when measurement and control are performed over a packet-dropping link \cite{gupta2007optimal}.
Separating estimation from control is practically convenient, and it is a commonly used paradigm in practice, even when it is not guaranteed to be optimal.

We seek to develop a fundamental understanding of the trade-offs which arise when applying the  separation principle to nonlinear and partially observed settings.
We show that in the presence of bilinear observations, the optimal control law does not obey the separation principle, nor is it linear.
Our results highlight 
the potential for tension between control and measurement, especially when the observations are less noisy than the state evolution.  
In particular, we show that the minimum of the LQG cost may maximize the state estimation error by causing a loss of observability.
Such input-dependent observability is a key challenge which arises in nonlinear and partially observed settings.

Classic optimal control problems for linear dynamics and measurements are well studied.
Perhaps the most well known is the $\mathcal H_2$ criterion, which results in a linear controller that satisfies the separation principle \cite{atassi1999separation, wonham1968separation}.
The robust $\mathcal H_\infty$ criterion also leads to a linear controller, and it asymptotically follows the separation principle~\cite{doyle1988state}.
Linear controllers are optimal in other settings, such as the state feedback min/max problem proposed by~\cite{gurpegui2023minimax}.
Even when not optimal, linear controllers are widely used.
Modern approaches based on
online learning~\cite{hazan2022introduction, simchowitz2019learning} and system level synthesis~\cite{anderson2019system}
restrict controllers to be linear, mainly out of convenience, but do not enforce or utilize the separation principle.

Bilinear optimal control is generally intractable \cite{pardalos2010optimization}, thus a natural focus for control synthesis is stabilization. To achieve stabilization, performing state estimation is a key step. For fixed input sequences, optimal state estimation for bilinear systems
follows the same principles as linear time varying systems.
As a result, when the noise is Gaussian, Kalman filtering is statistically optimal~\cite{kalman1960new, bertsekas2012dynamic}.
However, due to the input-dependent dynamics matrices, the Kalman filter is \emph{nonlinear} in the input and output variables.
In fact, unlike for linear systems, different choices of inputs may lead to better or worse filter performance.
The problem of selecting inputs to minimize estimation error is traditionally studied in active learning or experiment design~\cite{goodwin1977dynamic, pukelsheim2006optimal, bombois2011optimal}; however these settings treat the underlying state as static.
Considering how best to measure the state of a dynamical systems is also studied in sensor design~\cite{tanaka2016semidefinite}.
However, work in this area often focuses on discrete and combinatorial decision spaces, and do not model the ``observer effect'' of measurement on the state.

The relationship between estimation and control also plays an important role in optimal decentralized control. The famous Witsenhausen counterexample \cite{witsenhausen68} constructed a simple two-stage LQG system where two (scalar valued) control decisions are made given decentralized information, and showed that for this system, linear controllers are not optimal. This simple counterexample has played a crucial role in  understanding the challenges of decentralized control.

{\emph{Contributions:}} We make the following contributions towards addressing the problem of optimal control from bilinear observations: 
(i) we show that the separation principle does not hold when the observation model is bilinear, and (ii) provide a negative result, stating that the optimal controller is not affine in the estimated state. 
(iii) Under certain settings, we show that the optimal LQG controller can result in significantly worse performance, as it locally maximizes the cost function instead of minimizing. 
(iv) We then derive an analytical expression for the optimal nonlinear controller in these settings. (v) We also introduce a notion of input dependent observability and derive conditions under which the Kalman filtering covariance remains bounded.

\section{Preliminaries and Problem Setting}
{\emph{Notations:}} We use boldface uppercase (lowercase) letters to denote matrices (vectors). 
We use normal letters (uppercase \& lowercase) to denote scalars. 
For a vector $\vu$, $(\vu)_i$ denotes its $i$-th element. 
For a matrix $\vX$, $\norm{\vX}$ denotes its spectral norm, whereas, $\tr(\vX)$ denotes the trace of $\vX$.

\subsection{Bilinear Observation Optimal Control Problem}
In this paper, we consider the problem of optimally controlling a partially observed linear dynamical system~(LDS) from bilinear observations~(BO), with the objective of minimizing quadratic costs in states and actions. Specifically, we consider the following state-space representation: for all $t \geq 0$
\begin{equation}
\begin{aligned} \label{eqn: dynamics update}
    	\xb_{t+1} &= \Ab \xb_t + \Bb\ub_t + \wb_{t}, \\
\vy_t &= \big(\vC_0 + \sum_{k=1}^p(\vu_t)_k\vC_k\big) \xb_t + \vz_t, 
\end{aligned}
\end{equation}
where $\xb_t \in \R^n$ is the state, $\ub_t \in \R^p$ is the input, $\vy_t \in \R^m$ is the output, $\wb_t \in \R^n$ is the process noise, and $\vz_t \in \R^m$ is the measurement noise at time $t$. The matrices $\vA \in \R^{n \times n}$, and $\vB \in \R^{n \times p}$ govern the state evolution, whereas, $\vC_0, \vC_1, \dots, \vC_p \in \R^{m \times n}$ govern the bilinear observation model. For notational convenience, we define the input dependent observation matrix $\vC(\vu_t) := \vC_0 + \sum_{k=1}^p(\vu_t)_k\vC_k$. In the remainder of the paper, we use LDS-BO$(\vA, \vB,\vC_0, \vC_1, \dots, \vC_p)$ to refer to the state-space representation~\eqref{eqn: dynamics update}. With these definitions, we want to analyze the following finite-horizon optimal control problem:
\begin{equation}
\begin{aligned}\label{eqn:bilinear LQG}
    \inf_{\vu_{0}, \dots, \vu_{T-1}} &\E \left[ \vx_T^\T \vQ_T \vx_T + \sum_{t=0}^{T-1} \left( \vx_t^\T \vQ \vx_t + \vu_t^\T \vR \vu_t \right) \right],  \\
	&\text{s.t.} \quad \text{LDS-BO}(\vA, \vB,\vC_0, \vC_1, \dots, \vC_p),
\end{aligned}
\end{equation}
where $\vQ_T, \vQ, \vR$ are the positive semi-definite cost matrices which are assumed to be known. In order to solve the optimal control problem~\eqref{eqn:bilinear LQG}, we assume that the dynamic matrices $\vA,\vB, \vC_0,\vC_1, \dots, \vC_p$ are also known a priori. Moreover, the initial state, noise processes and cost matrices satisfy the following properties:
 \begin{assumption}\label{assump noise, initial, cost} We have: (i) Gaussian noise processes, $\{\vw_t\}_{t=0}^{T-1} \distas \Ncal(0, \vSigma_w)$, $\{\vz_t\}_{t=0}^{T-1} \distas \Ncal(0, \vSigma_z)$, (ii) initial state distributed as $\vx_0 \distas \Ncal(\vxhat_0, \vSigma_0)$, and (iii) the positive definite cost matrices, i.e., $\vQ_T, \vQ, \vR \succ 0$.    
 \end{assumption}
Note that Assumption~\ref{assump noise, initial, cost} is standard in the Linear Quadratic Gaussian~(LQG) control problem~\cite{aastrom2012introduction,mania2019certainty}, which is a special case of \eqref{eqn:bilinear LQG} with $\{\vC_k\}_{k=1}^p = \mathbf{0}_{mn}$.
Since we don't have access to the state vector $\vx_t$ at any time $t$, our control policy will be a function of all the information available at time $t$.
\begin{definition}[Admissible Policy] \label{def:admissible policy}
    Starting from $\Ical_0 = \phi$, let $\Ical_t := \{\vu_0, \vu_1, \dots, \vu_{t-1} ,\vy_0, \vy_1, \dots, \vy_{t-1} \}$ be the information available at time $t$ with $\Ical_{t+1} = \{\Ical_t, \vu_t, \vy_t\}$. Then, a policy $\mu_k$ is admissible if $\vu_t = \mu_t(\Ical_t)$.
\end{definition}
Note that, unlike the standard LQG case, the information vector $\Ical_t$ does not include $\vy_t$. This is because in~\eqref{eqn: dynamics update} one must take the action $\vu_t$ to observe the output $\vy_t$. In the next section, we will show how this effects the state estimation $\E[\vx_t |\Ical_t]$ at time $t$. Moreover, we will also show that the state estimation $\E[\vx_t |\Ical_t]$, and the estimation error covariance $\E[(\vx_t - \E[\vx_t |\Ical_t])(\vx_t - \E[\vx_t |\Ical_t])^\T \big|\Ical_t]$ are the sufficient statistics to design the optimal feedback control policy.

\subsection{Input-Dependent State Estimation}\label{sec:input_dep_state_est}
In this section, we will show the optimality of Kalman Filtering Algorithm~\cite{kalman1960new} for predicting the state of the dynamical system~\eqref{eqn: dynamics update}. For the sake of completeness, we first present a variant of Kalman Filtering~(KF) Algorithm here:

\begin{mdframed}\label{def:Kalman Filtering}
  {\bf Bilinear Observation Kalman Filtering:} Let $\vxhat_{t|t-1} := \E [\vx_t | \Ical_t]$ be the estimated state at time $t$, and $\vSigma_{t|t- 1} := \E \left[(\vx_t - \E [\vx_t | \Ical_t])(\vx_t - \E [\vx_t | \Ical_t])^\T\right]$ be the estimation error covariance. Then, starting from $\vxhat_{0|- 1} = \vxhat_0$, and $\vSigma_{0|- 1} = \vSigma_0$, the estimated state, and the error covariance follows the recursion:
\begin{equation*}
\begin{aligned}
    \vxhat_{t+1|t} &= \vA \vxhat_{t|t-1} + \vB \vu_t - \vL(\vu_t)\left(\vy_t - \vC(\vu_t) \vxhat_{t|t-1} \right), \\
    \vSigma_{t+ 1|t} &= \vA \vSigma_{t|t- 1} \vA^\T  + \vL(\vu_t)\vC(\vu_t) \vSigma_{t|t- 1} \vA^\T + \vSigma_w,\\
    \text{where,} \quad \vL(\vu_t) &:= -  \vA \vSigma_{t|t- 1} \vC(\vu_t)^\T \left( \vC(\vu_t) \vSigma_{t|t- 1} \vC(\vu_t)^\T + \vSigma_z \right)^{-1}
\end{aligned}
\end{equation*}
is the input-dependent Kalman gain at time $t$.
\end{mdframed}
Please refer to \cite[Eqs. (E.39) -- (E.42)]{bertsekas2012dynamic} for the derivation of Kalman filtering equations above. 
Importantly, unlike Kalman filtering from linear measurements, both the Kalman gain $\vL(\vu_t)$, and the error covariance $\vSigma_{t+ 1|t}$ depend on the control inputs up to time $t$. 
Hence, the separation principle~(Def.~\ref{def:Separation principle}) is not valid when controlling the partially observed LDS from bilinear  observations~\eqref{eqn: dynamics update}. 
However, it is easy to show that, the Kalman filtering gives the optimal state estimation of \eqref{eqn: dynamics update}, under Assumption~\ref{assump noise, initial, cost}.

\begin{lemma}[Optimality of KF] \label{lemma:optimality of KF} Under Assumption~\ref{assump noise, initial, cost}(i),(ii), the Kalman filtering algorithm gives:
\begin{enumerate}[label=(\textbf{\alph*}),leftmargin=*,noitemsep,topsep=0pt]
    \item optimal state estimation, that is, $\vxhat_{t|t- 1} = \E [\vx_t | \Ical_t]$.
    \item posterior distribution of the state $\vx_t$ given the information vector $\Ical_t$. Specifically, it gives $\vx_t | \Ical_t \sim \Ncal(\vxhat_{t|t- 1}, \vSigma_{t|t- 1})$.
\end{enumerate}
\end{lemma}
Lemma~\ref{lemma:optimality of KF} readily follows from~\cite{kalman1960new} by replacing $\vC_t$ with $\vC(\vu_t)$ which is a deterministic quantity conditioned on $\vu_t$. Moreover, since the posterior distribution of the state $\vx_t$ is Gaussian, using similar argument as~\cite[Section 4.3.1]{bertsekas2012dynamic}, we note that the values $\vxhat_{t|t- 1}, \vSigma_{t|t- 1}$ computed by the Kalman filter are sufficient statistics for any policy at time $t$.

\section{Main Results}
In this section we present our main results on the sub-optimality of the separation principle and linear controllers, and then derive the optimal nonlinear controller in a simple setting. First we present a formal definition of the separation principle, adapted from~\cite{aastrom2012introduction}.
\begin{definition}[Separation principle]\label{def:Separation principle}
    The optimal control from partial observations can be generated by first using Kalman filtering to find the optimal state estimate, then using linear quadratic regulator~(LQR) to find the optimal state feedback controller based on the optimal state estimate. 
\end{definition}
For the sake of completeness, in the following we present the sub-optimal LQG controller~\cite{bertsekas2012dynamic} for the system~\eqref{eqn: dynamics update}: at every time $0 \leq t \leq T-1$, the LQG controller is given by the sub-optimal policy:
\begin{equation}
\begin{aligned}\label{eqn:LQG Policy}
    \vu_t^{\rm LQG} &= \vL_t \vxhat_{t|t- 1} \quad \text{where,}\quad \vL_t &= -(\vB^\T \vK_{t+1} \vB + \vR)^{-1}\vB^\T \vK_{t+1} \vA,
\end{aligned}
\end{equation}
where the matrices $\vK_t$, starting from $\vK_T = \vQ_T$, are given recursively by the Riccati equation,
\begin{equation}
\begin{aligned}\label{eqn:LQG Riccati}
     \vP_t &= \vA^\T \vK_{t+1} \vB (\vB^\T \vK_{t+1} \vB + \vR)^{-1}\vB^\T \vK_{t+1} \vA, \\
     \vK_t &= \vA^\T \vK_{t+1} \vA - \vP_t + \vQ.
\end{aligned}
\end{equation}
We are now ready to present our main results.

\subsection{Sub-optimality of the Separation Principle} 
In this section, we will analyze the dynamic programming~(DP) algorithm to solve the optimal control problem~\eqref{eqn:bilinear LQG} as follows: Using the same proof technique as \cite[Section 4.2]{bertsekas2012dynamic}, the optimal policy for the last stage is given by $\vu_{T-1}^\star = \vu_{T-1}^{\rm LQG}$, which is the LQG controller given by~\eqref{eqn:LQG Policy}. However, the LQG controller is no longer optimal for $t \leq T-2$. To show this, we use the DP algorithm~\cite[Eq. (4.4)]{bertsekas2012dynamic} along-with the Woodbury matrix identity~\cite{woodbury1950inverting} to obtain
\begin{equation}
\begin{aligned} \label{eqn:u_Tm2_Bellman}
    \vu_{T-2}^\star = \arg \min_{\vu_{T-2}}  f(\vu_{T-2}) &:= f_{\rm LQG}(\vu_{T-2}) + g(\vu_{T-2}),\\
      \text{where} \quad f_{\rm LQG}(\vu_{T-2}) &:= \vu_{T-2}^\T \Acal \vu_{T-2} + 2 \vxhat_{T-2|T-3}^\T\Bcal^\T \vu_{T-2}, \\
      \text{and} \;\;\quad\quad\quad g(\vu_{T-2}) &:=  \tr \big(\Gcal \big(\vSigma_{T-2|T-3}^{-1} + \vC(\vu_{T-2})^\T \vSigma_z^{-1}\vC(\vu_{T-2})\big)^{-1} \big),
\end{aligned}
\end{equation}
where we define $\Acal:=\vB^\T \vK_{T-1}\vB + \vR$, $\Bcal:= \vB^\T \vK_{T-1} \vA$, and $\Gcal:=\vA^\T\vP_{T-1}\vA$ for notational convenience. 
Note that the second term $g(\vu_{T-2})$ depends on control input $\vu_{T-2}$ via the input dependent observation matrix $\vC(\vu_{T-2})$, and it corresponds to the estimation error covariance at time $T-1$. Hence, the control input at time $T-2$ directly affects the estimation error covariance at time $T-1$, violating the separation principle. 
Our next result shows that, because of the additional term $g(\vu_{T-2})$, the optimal control policy is not affine in the estimated state.

\begin{theorem}[Counter example] \label{thrm:negative result}
Consider the optimal control problem~\eqref{eqn:bilinear LQG}.  Suppose Assumption~\ref{assump noise, initial, cost} holds, and $T\geq 2$. Let $\{u^\star_k = \mu_k^\star(\Ical_k)\}_{k=0}^{T-1}$ be the optimal control policy obtained by minimizing~\eqref{eqn:bilinear LQG}. Then, for $k \leq T-2$, the optimal control policy $\mu_k^\star(\Ical_k)$ is not affine in the estimated state $\E[x_k | \Ical_k]$ in general.
\end{theorem}
The proof of Theorem~\ref{thrm:negative result} is presented in Section~\ref{proof of thrm:negative result}. Theorem~\ref{thrm:negative result} suggests that the optimal controller can be nonlinear under bilinear observations. However, finding an analytical expression for the optimal controller is challenging in general because of: (i) non-convexity, (ii) nonlinearity, and (iii) the existence of multiple maxima, minima, or saddle points. In the simple case of $T=2$, we find the optimal nonlinear controllers as follows.

\begin{theorem}[Nonlinear controller] \label{thrm:nonlinear special case}
    Consider the optimal control problem~\eqref{eqn:bilinear LQG}. 
    Suppose, $n=m=p=1$, and $T=2$. Let $P_{1}={A^2 B^2 Q_T^2}/(B^2 Q_T+R)$ and $K_{1} = {A^2 R Q_T}/(B^2 Q_T + R) + Q$.
    Suppose Assumption~\ref{assump noise, initial, cost} holds, and $\Sigma_z \leq C_1^2 {A^2 \Sigma_0^2 P_{T-1}}/{(B^2 K_{T-1} + R)}$. 
    Let $u^\star_0$, and $u^\star_1$  be the optimal control actions  obtained by minimizing~\eqref{eqn:bilinear LQG}. Let $u_0^{\rm LQG}$ be the optimal LQG control (obtained when $C_1 = 0$). Suppose $C_0 = -u_0^{\rm LQG} \times C_1$. Then: 
    \begin{enumerate}[label=(\textbf{\alph*}),leftmargin=*,noitemsep,topsep=0pt]
    \item $u_0^{\rm LQG}$ locally maximizes the cost in \eqref{eqn:bilinear LQG}.
    \item the optimal control action $u_0^\star$ is nonlinear in $\xhat_0$ and is not unique, that is, 
    \begin{align}
        u_0^\star &= - \frac{\beta \xhat_0}{\alpha} \left( 1 \pm \frac{1}{C_0} \sqrt{-\kappa + \frac{\alpha C_0}{\beta \xhat_0} \sqrt{\frac{\gamma}{\alpha}}} \right), \\
        u_1^\star &= - \left(\frac{A Q_T B}{B^2 Q_T + R}\right)  \xhat_{1|0},
    \end{align}
\end{enumerate}
where $\alpha := B^2 K_1 + R$, $\beta := B K_1 A$, $\gamma := \Sigma_z A^2 P_1$, and $\kappa := \Sigma_z\Sigma_0^{-1}$.
\end{theorem}
The proof of Theorem~\ref{thrm:nonlinear special case} is given in Section~\ref{Proof of thrm:nonlinear case}. 
The setting of Theorem~\ref{thrm:nonlinear special case} (supposing $C_0 = -u_0^{\rm LQG} \times C_1$) means that the LQG control input exactly cancels the state observation in~\eqref{eqn: dynamics update}, and thus maximizes the estimation error covariance.
This suggests that the landscape of the cost-to-go function depends on the distance between the LQG controller and the maximizer of the estimation error covariance. 
We investigate this further numerically (Figure~\ref{figure1}),
and see that
when this distance is large, the LQG controller is also the optimal controller. However, when this distance is small, the LQG controller becomes sub-optimal, and the optimal controllers are nonlinear is the estimated state.

\subsection{Sufficient Conditions for Uniformly Bounded Cost}
In this section, we derive conditions under which the estimation error covariance stays bounded. Since our observation matrices are input dependent, we first introduce the notion of uniform observability for bilinear observation systems~\ref{eqn: dynamics update}. 
\begin{definition}[Complete observability]\label{def:observability cond}
    The system \eqref{eqn: dynamics update} is \emph{uniformly completely observable} for input sequence $\{\vu_t\}_{t\geq 0}$ if there exists a $\delta>0$ such that, for all $\ell\geq 0$,
    \[\vO_\ell := \sum_{k=0}^{n-1} (\vA^k)^\T \vC(\vu_{\ell+k})^\T \vC(\vu_{\ell+k})\vA^k \succ \delta \vI_n.\]
\end{definition}
We call $\vO_\ell$ observability grammian related to~\eqref{eqn: dynamics update}. Note that unlike linear observation case, here  $\vO_\ell$ is time varying because of the input dependent observation matrix. Our next result~(adapted from~\cite{gurumoorthy2017rank}) establishes a useful connection between the KF error covariance and the observability condition introduced in Definition~\ref{def:observability cond}.

\begin{lemma}[Bounded error covariance]\label{lemma:bounded error covariance}
    Suppose that the system \eqref{eqn: dynamics update} is uniformly completely observable for input sequence $\{\vu_t\}_{t\geq 0}$. Then the Kalman filtering error covariance matrix remains bounded for all time, i.e. there exists a constant $c$ such that $\|\vSigma_{t|t-1}\|\leq c$ for all $t\geq 0$.
\end{lemma}
Lemma~\ref{lemma:bounded error covariance} follows readily from Lemma 3.2 in \cite{gurumoorthy2017rank}, noting that the positive definite condition in Definition~\ref{def:observability cond} implies the nonzero determinant condition in Definition 3.1 of \cite{gurumoorthy2017rank}.

The separation principle can lead to degenerate policies which fail to adequately estimate the state.
Consider, for example, settings with $\vC_0=\mathbf{0}_{mn}$.
When $\vxhat_0 = \E[\vx_0] = 0$, the system~\eqref{eqn: dynamics update} will be non-observable, and it is impossible to find a stabilizing control policy. This is because when $\vC_0=\mathbf{0}_{mn}$, $\vxhat_0 = 0$ implies $\vu_0 =0$, and $\vC(\vu_0) = \sum_{k=1}^p (\vu_0)_k \vC_k =\mathbf{0}_{mn}$. This will result in the Kalman gain $\vL(\vu_t)$ being zero for all $t \geq 0$. Hence, the estimated state is $\vxhat_t = 0$ for all $t\geq0$. 
Such failures may be preventable, e.g. by randomizing the initial state estimate.
While it is difficult to rule out this problem in general, for certain classes of observation models it is avoidable.
The following proposition~(proved in Section~\ref{proof of proposition perp}) provides sufficient conditions for ensuring observability regardless of the control policy.

\begin{proposition}[Sufficient condition] \label{prop:perpendicular}
    Let $\vC_0^\perp$ be the projection of $\vC_0$ onto the orthogonal complement of the span of $\vC_1,...,\vC_p$. Then if $(\vA,\vC_0^\perp)$ is observable, 
    for any choice of control inputs $\{\vu_t\}_{t\geq 0}$, the Kalman filtering error covariance matrix remains bounded for all $t{\geq0}$.
\end{proposition}

\section{Proofs of Main Results}

\subsection{Proof of Theorem~\ref{thrm:negative result}}\label{proof of thrm:negative result}
\begin{proof}
    Our proof is based on the analysis of the dynamic programming~(DP) algorithm to solve \eqref{eqn:bilinear LQG}. Fix $n=m=p=1$, and consider the following optimization problem,
\begin{equation}
\begin{aligned}\label{eqn:bilinear LQG scalar}
    \inf_{u_{0}, \dots, u_{T-1}} &\E \left[ Q_T x_T^2 + \sum_{t=0}^{T-1} \left( Q x_t^2 + R u_t^2 \right) \right],  \\
	&\text{s.t.} \quad \text{LDS-BO}(A, B, C_0, C_1).
\end{aligned}
\end{equation}
    To begin, we use~\cite[Eq. (4.4)]{bertsekas2012dynamic} to find the optimal controller $u_{T-1}^\star = \mu_{T-1}^\star(\Ical_{T-1})$ as follows,
\begin{align}
    u^\star_{T-1} &= - \left(\frac{A Q_T B}{B^2 Q_T + R}\right) \E \left[x_{T-1}  | \Ical_{T-1}\right], 
\end{align}
and the corresponding optimal cost for the last two time steps is given by
\begin{align}
    J_{T-1}^\star(\Ical_{T-1}) &= K_{T-1} {\E} \left[ x_{T-1}^2  | \Ical_{T-1}\right] +  Q_T \Sigma_w
    + P_{T-1}  {\E} \left[ \left(x_{T-1} - \E \left[x_{T-1}   | \Ical_{T-1}\right] \right)^2  | \Ical_{T-1}\right] , \nn
\end{align}
where $P_{T-1}={A^2 B^2 Q_T^2}/(B^2 Q_T+R)$ and $K_{T-1}= A^2Q_T - P_{T-1} + Q  = {A^2 R Q_T}/(B^2 Q_T + R) + Q$. Then using the DP algorithm~\cite[Eq. (4.8)]{bertsekas2012dynamic}, we have
\begin{align}
    J_{T-2}^\star(\Ical_{T-2}) 
    &= \min_{u_{T-2}} \bigg\{ {\E} \bigg[ Q x_{T-2}^2  + R u_{T-2}^2 +  Q_T \Sigma_w 
    + P_{T-1}  {\E} \big[ \left(x_{T-1} - \E \left[x_{T-1}   | \Ical_{T-1}\right] \right)^2  | \Ical_{T-1}\big]   \nn \\
    &+ K_{T-1} {\E} \left[ x_{T-1}^2  | \Ical_{T-1}\right] \big| \Ical_{T-2}, u_{T-2}\bigg] \bigg\}. \label{eqn:J_Tm2 expression}
\end{align}
Unlike the LQG controller, here we will show that the estimation error also depends on the control input, hence, violating the separation principle in Definition~\ref{def:Separation principle}. Specifically, using the Kalman filtering equations in Section~\ref{def:Kalman Filtering}, we obtain,
\begin{align}
    &{\E}\left[ {\E} \left[ \left(x_{T-1} - \E \left[x_{T-1} | \Ical_{T-1}\right] \right)^2  | \Ical_{T-1} \right] \big| \Ical_{T-2}, u_{T-2} \right] \nn \\
    &\eqsym{i}  {\E} \left[ \left(x_{T-1} - \xhat_{T-1|T-2} \right)^2   \big| \Ical_{T-2}, u_{T-2} \right], \nn \\
    & =:  \Sigma_{T-1 | T-2} \big| \Ical_{T-2}, u_{T-2} , \nn \\
    &\eqsym{ii}  A^2\Sigma_{T-2 | T-3}  -  \frac{A^2 (C_0 + C_1 u_{T-2})^2 \Sigma_{T-2|T-3}^2}{(C_0 + C_1 u_{T-2})^2 \Sigma_{T-2|T-3} + \Sigma_z} + \Sigma_w , \nn \\ 
    &=  \frac{ \Sigma_z A^2 P_{T-1}}{(C_0 + C_1 u_{T-2})^2  + \Sigma_z\Sigma_{T-2|T-3}^{-1}} + \Sigma_w, \label{eqn:error variance u_Tm2}
\end{align}
where we use the tower property along with the fact that $\{\Ical_{T-2}, u_{T-2}\} \subset \Ical_{T-1}$ to obtain (i), and (ii) follows from the Kalman filtering update. Combining \eqref{eqn:J_Tm2 expression}, and \eqref{eqn:error variance u_Tm2}, the optimal controller $u_{T-2}^\star = \mu_{T-2}^\star(\Ical_{T-2})$ is obtained by solving the following minimization problem,
\begin{equation}
\begin{aligned} \label{eqn:f of u_Tm2 expression}
    u_{T-2}^\star = \arg \min_{u} f_{T-2}(u):= \alpha u^2 + 2 \beta \xhat u +
    \frac{\gamma}{(C_0 + C_1 u)^2  + \kappa},
\end{aligned}
\end{equation}
where we define the following for notational convenience,
\begin{equation}
    \begin{aligned} \label{eqn:alpha beta gamma kappa}
    \alpha &:= B^2 K_{T-1} + R \qquad \qquad \beta := B K_{T-1} A \\
    \gamma &:= \Sigma_z A^2 P_{T-1} \qquad \qquad~~
    \kappa := \Sigma_z\Sigma_{T-2|T-3}^{-1}
\end{aligned}
\end{equation}
Moreover, we set $\xhat := \xhat_{T-2|T-3}$, and $u := u_{T-2}$ for notational convenience. In the following, we will use the proof by contradiction to show that $u_{T-2}^\star$ is not affine in $\xhat_{T-2|T-3}$. To begin, we do a change of variable by setting $\ubar = C_0 + C_1 u  \implies u = (\ubar-C_0)/C_1$. Then, the minimization problem~\eqref{eqn:f of u_Tm2 expression} becomes,
\begin{equation}
\begin{aligned} \label{eqn:f bar of u_Tm2 expression}
    \ubar_{T-2}^\star = \arg \min_{\ubar} \fbar_{T-2}(\ubar)&:= \alpha \left(\frac{\ubar - C_0}{C_1}\right)^2 
    + 2 \beta \xhat \left(\frac{\ubar - C_0}{C_1}\right) +  \frac{\gamma}{\ubar^2  + \kappa},
\end{aligned}
\end{equation}
Setting the derivative of $\fbar_{T-2}(\ubar)$ equal to zero, we obtain the following $5$-th order polynomial equation with one of the roots corresponding to $\ubar_{T-2}^\star$,
\begin{align}
&\alpha \ubar^5 + (C_1\beta \xhat - C_0 \alpha) \ubar^4 + 2 \kappa \alpha \ubar^3 \nn  \\
    &+ 2\kappa (C_1\beta \xhat - C_0 \alpha) \ubar^2
    + (\alpha \kappa^2 - \gamma C_1^2) \ubar
    + \kappa^2 (C_1\beta \xhat - C_0 \alpha) = 0. \label{eqn:Gradient condition}
\end{align}
Suppose there exists some $0 \neq \Lbar \in \R $, and $\dbar \in \R$ such that the optimal controller can be represented as $\ubar_{T-2}^\star = \Lbar\xhat + \dbar$. Then, from \eqref{eqn:Gradient condition}, we have
\begin{align}
     & \alpha (\Lbar\xhat + \dbar)^5 + (C_1\beta \xhat - C_0 \alpha) (\Lbar\xhat + \dbar)^4 + 2 \kappa \alpha (\Lbar\xhat + \dbar)^3  \nn \\
     &+ 2\kappa (C_1\beta \xhat - C_0 \alpha) (\Lbar\xhat + \dbar)^2 + (\alpha \kappa^2 - \gamma C_1^2) (\Lbar\xhat + \dbar) 
     + \kappa^2 (C_1\beta \xhat - C_0 \alpha) = 0.
\label{eqn:affine_contril_cond}
\end{align}
Since this holds for all $\xhat \in \R$, \eqref{eqn:affine_contril_cond} implies that we must have $\lambda_i = 0$ for all $i = 0,1, \dots, 5$, where $\lambda_i$ is the coefficient of $\xhat^i$ in \eqref{eqn:affine_contril_cond}. Utilizing this information, for $i=5$, we have
\begin{align}
    \lambda_5 = \alpha \Lbar^5 + C_1 \beta \Lbar^4  = 0 \implies \Lbar = -\frac{C_1 \beta}{\alpha}. \label{eqn:Lbar Lambda_5}
\end{align}
Similarly, for $i=4$, using \eqref{eqn:Lbar Lambda_5}, we have
\begin{align}
    \lambda_4 &= 5 \alpha \dbar \Lbar^4 - C_0 \alpha \Lbar^4 + 4 C_1 \beta \dbar \Lbar^3, \nn \\
    &= \Lbar^3 \left( C_0 C_1 \beta  -  C_1 \beta \dbar\right) = 0 \implies \dbar = C_0 \label{eqn:dbar Lambda_4}
\end{align}
Finally, for $i=0$, using \eqref{eqn:dbar Lambda_4}, we have
\begin{align}
     \lambda_0 &= \alpha \dbar^5 - \alpha C_0 \dbar^4 + 2 \kappa \alpha \dbar^3 -  2 \kappa \alpha C_0\dbar^2 + \kappa^2 \alpha \dbar - \kappa^2 \alpha C_0 - C_1^2 \dbar \gamma \nn \\
     &   = - C_1^2 C_0 \gamma \neq 0, \label{eqn:lambda_0 contradiction}  
\end{align}
under the condition that $A, B , Q_T, C_0, C_1$, and $\Sigma_z$ are non-zero. Alternately, we can also show the contradiction by noticing that,
\begin{align}
    \lambda_1 &= 5 \alpha \Lbar \dbar^4 + C_1 \beta\dbar^4 - 4C_0 \alpha \Lbar \dbar^3 + 6\kappa\alpha \Lbar \dbar^2 + 2 \kappa C_1 \beta \dbar^2 \nn \\
    &- 4 \kappa C_0 \alpha \Lbar \dbar + (\alpha \kappa^2 - \gamma C_1^2) \Lbar + \kappa^2 C_1 \beta, \nn  \\
&\eqsym{i} (\alpha \Lbar + C_1 \beta) C_0^4 + 2\kappa (\alpha \Lbar  + C_1 \beta) C_0^2 
    + \kappa^2 (\alpha \Lbar + C_1 \beta) - \gamma C_1^2 \Lbar, \nn \\
    & \eqsym{ii} \frac{\gamma C_1^3 \beta}{\alpha} \neq 0, \label{eqn:lambda_1 contradiction}  
\end{align}
under the condition that $A, B , Q_T, Q, R, C_1$, and $\Sigma_z$ are non-zero. This contradicts our claim that $\lambda_i = 0$ for all $i = 0,1, \dots, 5$, and completes the proof. 
\end{proof}

\subsection{Proof of Theorem~\ref{thrm:nonlinear special case}} \label{Proof of thrm:nonlinear case}
\begin{proof}
To begin, recall from the proof of Theorem~\ref{thrm:negative result}, that the optimal controller $u_{T-2}^\star = \mu_{T-2}^\star(\Ical_{T-2})$ is obtained by solving the minimization problem~\eqref{eqn:f of u_Tm2 expression}. When $T=2$, we have $u_0^\star = \arg \min_{u} f_0(u)$, where
\begin{align}
    f_0(u) &= \alpha u^2 + 2 \beta \xhat_0 u +  \frac{\gamma}{(C_0 + C_1 u)^2  + \kappa}, \label{eqn: f of zero expression} 
\end{align}
where $\alpha, \beta, \gamma$, and $\kappa$ are as defined in \eqref{eqn:alpha beta gamma kappa}. Similar to the proof of Theorem~\ref{thrm:negative result}, we do a change of variable by setting $\ubar = C_0 + C_1 u  \implies u = (\ubar-C_0)/C_1$. Then, the minimization problem  becomes, $\ubar_0^\star = \arg \min_{\ubar} \fbar_0(u)$, where
\begin{align*}
    \fbar_0(\ubar) &= \alpha \left(\frac{\ubar - C_0}{C_1}\right)^2 + 2 \beta \xhat_0 \left(\frac{\ubar - C_0}{C_1}\right) +  \frac{\gamma}{\ubar^2  + \kappa}, \\
    \fbar_0'(\ubar) &= \frac{2\alpha}{C_1} \left(\frac{\ubar - C_0}{C_1}\right)  +  \frac{2 \beta \xhat_0}{C_1} -  \frac{2\gamma \ubar}{(\ubar^2+\kappa)^2}, \\
    \fbar_0''(\ubar) &=  \frac{2\alpha}{C_1^2} +  \frac{2\gamma \left( 3 \ubar^2 - \kappa\right)}{(\ubar^2  + \kappa)^3}. 
\end{align*}
Setting the gradient $\fbar_0'(\ubar) = 0$ equal to zero, we have
\begin{align}
     (\ubar^2+\kappa)^2(\alpha \ubar +  C_1\beta \xhat_0 - C_0 \alpha) - \gamma C_1^2 \ubar = 0. \label{eqn:gradient equal zero thrm 2}
\end{align}
Using the assumption that $-\frac{C_0}{C_1} = -\frac{\beta \xhat_0}{\alpha} =: u_0^{LQG}$, we get $C_1\beta \xhat_0 - C_0 \alpha = 0$. This further implies that the optimal $\ubar_0^\star$ is one of the roots of the following polynomial,
\begin{align}
    \left(\alpha (\ubar^2+\kappa)^2  - \gamma C_1^2 \right) \ubar = 0 \label{eqn:root polynomial thrm 2} 
\end{align}
which are found as follows:

    {\bf (a)} $\ubar_0^\star = 0$ or $u_0^\star = -\frac{C_0}{C_1} =  -\frac{\beta \xhat_0}{\alpha} = u_0^{LQG}$ which can be a maxima, minima or a saddle point depending on the value of, 
    \begin{align}
        \fbar_0''(\ubar_0^\star) = \frac{2\alpha}{C_1^2} -  \frac{2\gamma}{\kappa^2} =  2\frac{\alpha \kappa^2 - \gamma C_1^2}{C_1^2 \kappa^2}.
    \end{align}
    Specifically; {\bf (i)} when $\alpha \kappa^2 > \gamma C_1^2$, then $\fbar_0''(\ubar_0^\star) >0$, hence, $\ubar_0^\star = 0$ is the local minima of $\fbar_0(\ubar)$; {\bf (ii)} when $\alpha \kappa^2 < \gamma C_1^2$, then $\fbar_0''(\ubar_0^\star) <0$, hence, $\ubar_0^\star = 0$ is the local maxima of $\fbar_0(\ubar)$; {\bf (iii)} when $\alpha \kappa^2 = \gamma C_1^2$, then $\fbar_0''(\ubar_0^\star) =0$. However, we show that $\ubar_0^\star = 0$ is the local minima of $\fbar_0(\ubar)$ as follows: For any $\delta \in \R$, we have,
        \begin{align}
            \fbar_0'(\delta) &=  \frac{2\alpha}{C_1^2} \left(\delta - C_0\right)  +  \frac{2 \beta \xhat_0}{C_1} -  \frac{2\gamma \delta}{(\delta^2+\kappa)^2}
= 2 \delta \frac{\alpha \delta^4+ 2\alpha\delta^2 \kappa}{C_1^2(\delta^2+\kappa)^2}.
        \end{align}
        Hence, using the fact that $\alpha, \kappa >0$, we have $\lim_{\delta \to 0^-} \fbar_0'(\delta) < 0$ and $\lim_{\delta \to 0^+} \fbar_0'(\delta) > 0$. Hence, $\ubar_0^\star = 0$ is the local minima of $\fbar_0(\ubar)$.

    {\bf (b)} The roots of the following fourth order polynomial:
        \begin{align}
\alpha \ubar^4 + 2 \alpha\kappa \ubar^2 + \alpha \kappa^2 - \gamma C_1^2 = 0 
\implies \ubar_0^\star &=  \pm \sqrt{-\kappa \pm C_1 \sqrt{\frac{\gamma}{\alpha}}}, \nn\\
             \implies u_0^\star &=  - \frac{\beta \xhat_0}{\alpha} \pm \frac{1}{C_1}\sqrt{-\kappa \pm C_1 \sqrt{\frac{\gamma}{\alpha}}},
        \end{align}
    which can be a maxima, minima or saddle point depending on the value of,
    \begin{align}
        \fbar_0''\left(\pm \sqrt{-\kappa + C_1 \sqrt{\frac{\gamma}{\alpha}}}\right) &= 8 \gamma \frac{ C_1\sqrt{\frac{\gamma}{\alpha}} -  \kappa  }{C_1^3 (\frac{\gamma}{\alpha})^{3/2}}.
    \end{align}
    Specifically; {\bf (i)} when $C_1\sqrt{\frac{\gamma}{\alpha}} -  \kappa < 0$ (equivalently, $\alpha \kappa^2 > \gamma C_1^2$), then $\fbar_0''\left(\pm \sqrt{-\kappa + C_1 \sqrt{\frac{\gamma}{\alpha}}}\right) <0$, however, $\ubar_0^\star  =  \pm i\sqrt{\kappa - C_1 \sqrt{\frac{\gamma}{\alpha}}}$ are the complex roots of \eqref{eqn:root polynomial thrm 2}. Hence, $u_0^\star = - \frac{\beta \xhat_0}{\alpha} \pm i\frac{1}{C_1}\sqrt{\kappa - C_1 \sqrt{\frac{\gamma}{\alpha}}}$ are the complex roots of $f_0'(u)=0$; 
    {\bf (ii)} when $C_1\sqrt{\frac{\gamma}{\alpha}} -  \kappa > 0$ (equivalently, $\alpha \kappa^2 < \gamma C_1^2$), then $\fbar_0''\left(\pm \sqrt{-\kappa + C_1 \sqrt{\frac{\gamma}{\alpha}}}\right) >0$, hence, $\ubar_0^\star = \pm \sqrt{-\kappa + C_1 \sqrt{\frac{\gamma}{\alpha}}} $ are the local minima $f(\ubar)$. Hence, $u_0^\star = - \frac{\beta \xhat_0}{\alpha} \pm \frac{1}{C_1}\sqrt{-\kappa + C_1 \sqrt{\frac{\gamma}{\alpha}}}$ are the local minima $f(u)$; 
    {\bf (iii)} when $C_1\sqrt{\frac{\gamma}{\alpha}} -  \kappa = 0$ (equivalently, $\alpha \kappa^2 = \gamma C_1^2$), then $\ubar_0^\star = 0 $ is the local minima $f(\ubar)$. Hence, $u_0^\star = - \frac{\beta \xhat_0}{\alpha}$ is the local minima $f(u)$.

    {\bf (c)} The remaining two roots will always be the complex roots of~\eqref{eqn:root polynomial thrm 2}, and are given by $\ubar_0^\star = \pm i\sqrt{\kappa + C_1 \sqrt{\frac{\gamma}{\alpha}}}$. Hence, $u_0^\star = - \frac{\beta \xhat_0}{\alpha} \pm i\frac{1}{C_1}\sqrt{\kappa + C_1 \sqrt{\frac{\gamma}{\alpha}}}$ are the complex roots of $f_0'(u)=0$. 
    
    Finally, the proof of Theorem~\ref{thrm:nonlinear special case} follows from {\bf (a)}(ii), and {\bf (b)}(ii). This completes the proof.

\end{proof}

\subsection{Proof of Proposition~\ref{prop:perpendicular}} \label{proof of proposition perp}
\begin{proof}
    Consider Frobenius inner product between matrices: $\langle \vU, \vV \rangle = \tr(\vU^\T\vV)$.
    By the fundamental theorem of linear algebra, we can decompose 
    $\vC_0= \vC_0^\perp  + \bar \vC_0$  where $\bar \vC_0$ is in the span of $\vC_1,...,\vC_p$.
    Define $\bar\vC(\vu_t) = \vC(\vu_t)-\vC_0^\perp$, then we have $\tr\big((\vC_0^\perp)^\T \bar\vC(\vu_t)\big)=0$. 
    Then the observability Gramian $\vO_\ell$ for any $\ell$ can be decomposed into $\vO = \vO_1+\vO_2+\vO_3$ where
    \begin{align}
        \vO_1 &= \sum_{i=0}^{n-1} (\vA^i)^\T (\vC_0^\perp)^\T \vC_0^\perp\vA^i \nn\\
        \vO_2 &= \sum_{i=0}^{n-1} (\vA^i)^\T \left[ (\vC_0^\perp)^\T \bar\vC(\vu_{\ell+i}) +\bar\vC(\vu_{\ell+i})^\T \vC_0^\perp \right]\vA^i  \nn \\
        \vO_3 &= \sum_{i=0}^{n-1} (\vA^i)^\T \bar\vC(\vu_{\ell+i})^\T \bar\vC(\vu_{\ell+i})\vA^i
    \end{align}
    The first term $\vO_1$ is the time invariant observability Grammian for $(\vA,\vC_0^\perp)$, which is positive definite by assumption. 
    For $\vO_2$, observe that $\vv^\T\vO_2\vv=0$ for all $\vv \in \R^n$ because the following inequality holds for any $\vA,\vB \in \R^{m\times n}$ and $\vv \in \R^n$:
    \begin{align}
        |\vv^\T \vA^\T \vB \vv| &= |\tr\left( \vA^\T \vB \vv \vv^\T\right)| \nn \\
        &= \sqrt{\tr\left(\vA^\T \vB \vv\vv^\T\right)^2} \leq \sqrt{\tr(\vA^\T \vB)\, \tr(\vv\vv^\T)}
    \end{align}
    Lastly, the matrix $\vO_3$ is clearly positive semidefinite.
    This shows that $\vO$ is positive definite, and thus the system is uniformly completly observable for any choice of inputs.
    The conclusion follows by Lemma~\ref{lemma:bounded error covariance}.
\end{proof}
\subsection{Derivation of Equation~\eqref{eqn:u_Tm2_Bellman}}
\begin{proof}    
To begin, we use~\cite[Eq. (4.4)]{bertsekas2012dynamic} to find the optimal control input at time $T-1$ as follows,
\begin{align}
    \vu^\star_{T-1} =\mu_{T-1}^\star(\Ical_{T-1}) = - (\vB^\T\vQ_T\vB + \vR)^{-1} \vB^\T \vQ_T \vA \E \left[\vx_{T-1}  | \Ical_{T-1}\right], 
\end{align}
and the corresponding optimal cost for the last two time steps is given by
\begin{align}
    J_{T-1}^\star(\Ical_{T-1}) &= \E \left[\vx_{T-1}^\T \vK_{T-1} \vx_{T-1} | \Ical_{T-1}\right] + \tr \left(\vQ_T \vSigma_w\right) \nn \\
    &+  \E \left[ \left(\vx_{T-1} - \E \left[\vx_{T-1}   | \Ical_{T-1}\right] \right)^\T \vP_{T-1} \left(\vx_{T-1} - \E \left[\vx_{T-1}   | \Ical_{T-1}\right] \right) \big| \Ical_{T-1}\right],
\end{align}
where the matrices $\vP_{T-1}$ and $\vK_{T-1}$ are given by
\begin{align}
    \vP_{T-1} &= \vA^\T \vQ_T \vB (\vR + \vB^\T \vQ_T \vB)^{-1} \vB^\T \vQ_T \vA,  \nn \\
    \vK_{T-1} &= \vA^\T \vQ_T \vA - \vA^\T \vQ_T \vB (\vR + \vB^\T \vQ_T \vB)^{-1} \vB^\T \vQ_T \vA + \vQ.
\end{align}
Then using the DP algorithm~\cite[Eq. (4.8)]{bertsekas2012dynamic}, the optimal cost for the last three time steps is, 
\begin{align}
    J_{T-2}^\star(\Ical_{T-2}) 
    &= \min_{\vu_{T-2}} \bigg\{  \E \bigg[ \vx_{T-2}^\T \vQ \vx_{T-2} + \vu_{T-2}^\T \vR \vu_{T-2} + \E \left[\vx_{T-1}^\T \vK_{T-1} \vx_{T-1} | \Ical_{T-1}\right] + \tr \left(\vQ_T \vSigma_w\right) \nn \\
    &+  \E \left[ \left(\vx_{T-1} - \E \left[\vx_{T-1} | \Ical_{T-1}\right] \right)^\T \vP_{T-1} \left(\vx_{T-1} - \E \left[\vx_{T-1} | \Ical_{T-1}\right] \right) \big| \Ical_{T-1} \right] \bigg| \Ical_{T-2}, \vu_{T-2}\bigg]\bigg\}. \label{eqn:J_Tm2 expression_vector}
\end{align}
Unlike the standard LQG case, the last term in \eqref{eqn:J_Tm2 expression_vector} depends on the input $\vu_{T-2}$, as follows,
\begin{align}
    &\E \left[ \E \left[ \left(\vx_{T-1} - \E \left[\vx_{T-1} | \Ical_{T-1}\right] \right)^\T \vP_{T-1} \left(\vx_{T-1} - \E \left[\vx_{T-1} | \Ical_{T-1}\right] \right) \big| \Ical_{T-1} \right] \bigg| \Ical_{T-2}, \vu_{T-2} \right] \nn \\
    &=\E \left[ \E \left[ \left(\vx_{T-1} - \vxhat_{T-1|T-2} \right)^\T \vP_{T-1} \left(\vx_{T-1} - \vxhat_{T-1|T-2} \right) \big| \Ical_{T-1} \right] \bigg| \Ical_{T-2}, \vu_{T-2} \right], \nn \\
    &=\E\left[\tr\left(\vP_{T-1}\vSigma_{T-1|T-2}\right)\big| \Ical_{T-2}, \vu_{T-2}\right]. \label{eqn:additional term u_Tm2}
\end{align}
To proceed, using the bilinear observation Kalman filtering in Section~\ref{def:Kalman Filtering}, we get the following,
\begin{align}
    \vSigma_{T-1|T-2} &= \vA \bigg( \vSigma_{T-2|T-3}- \vSigma_{T-2|T-3}\vC(\vu_{T-2})^\T \left(\vC(\vu_{T-2})\vSigma_{T-2|T-3}\vC(\vu_{T-2})^\T + \vSigma_z\right)^{-1} \nn \\
    &\quad~ \vC(\vu_{T-2})\vSigma_{T-2|T-3}\bigg) \vA^\T + \vSigma_w, \nn \\
    &= \vA \vSigma_{T-2|T-3}^{1/2} \bigg( \vI_n - \vSigma_{T-2|T-3}^{1/2}\vC(\vu_{T-2})^\T \left(\vC(\vu_{T-2})\vSigma_{T-2|T-3}\vC(\vu_{T-2})^\T + \vSigma_z\right)^{-1} \nn \\
    &\quad~ \vC(\vu_{T-2})\vSigma_{T-2|T-3}^{1/2}\bigg) \vSigma_{T-2|T-3}^{1/2}\vA^\T + \vSigma_w, \nn \\
    &\eqsym{i} \vA\vSigma_{T-2|T-3}^{1/2} \left(\vI_n + \vSigma_{T-2|T-3}^{1/2}\vC(\vu_{T-2})^\T \vSigma_z^{-1}\vC(\vu_{T-2})\vSigma_{T-2|T-3}^{1/2}\right)^{-1} \vSigma_{T-2|T-3}^{1/2}\vA^\T + \vSigma_w, \nn \\
    &= \vA \left(\vSigma_{T-2|T-3}^{-1} + \vC(\vu_{T-2})^\T \vSigma_z^{-1}\vC(\vu_{T-2})\right)^{-1}\vA^\T + \vSigma_w,
\end{align}
where we obtained (i) by using Woodbury matrix identity~\cite{woodbury1950inverting}. Combining this with \eqref{eqn:additional term u_Tm2}, we have
\begin{align}
    &\E \left[ \E \left[ \left(\vx_{T-1} - \E \left[\vx_{T-1} | \Ical_{T-1}\right] \right)^\T \vP_{T-1} \left(\vx_{T-1} - \E \left[\vx_{T-1} | \Ical_{T-1}\right] \right) \big| \Ical_{T-1} \right] \bigg| \Ical_{T-2}, \vu_{T-2} \right] \nn \\
& = \tr \left(\vP_{T-1}\vA \left(\vSigma_{T-2|T-3}^{-1} + \vC(\vu_{T-2})^\T \vSigma_z^{-1}\vC(\vu_{T-2})\right)^{-1} \vA^\T\right) + \tr \left(\vP_{T-1}\vSigma_w\right), \nn \\
    & = \tr \left(\vA^\T\vP_{T-1}\vA \left(\vSigma_{T-2|T-3}^{-1} + \vC(\vu_{T-2})^\T \vSigma_z^{-1}\vC(\vu_{T-2})\right)^{-1} \right) + \tr \left(\vP_{T-1}\vSigma_w\right), \label{eqn:additional term u_Tm2 final}
\end{align}
Finally, using \eqref{eqn:additional term u_Tm2 final} in \eqref{eqn:J_Tm2 expression_vector}, and after some algebraic manipulations, the optimal control input at time $T-2$ is given by the following non-convex minimization problem,
\begin{align}
    \vu_{T-2}^\star = \arg \min_{\vu_{T-2}} \bigg\{ &\vu_{T-2}^\T \Ac \vu_{T-2} + 2 \vxhat_{T-2|T-3}^\T\Bc^\T \vu_{T-2} \nn \\ 
    &+ \tr \left(\Gc \left(\vSigma_{T-2|T-3}^{-1} + \vC(\vu_{T-2})^\T \vSigma_z^{-1}\vC(\vu_{T-2})\right)^{-1} \right)\bigg\}
\end{align}
where we set $\Acal:=\vB^\T \vK_{T-1}\vB + \vR$, $\Bcal:= \vB^\T \vK_{T-1} \vA$, and $\Gcal:=\vA^\T\vP_{T-1}\vA$ for notational convenience. This completes the proof. 
\end{proof}

\section{Numerical Experiments}

\begin{figure*}[t]
\setlength{\abovecaptionskip}{0pt}
\setlength{\belowcaptionskip}{0pt}
\begin{center}
\begin{tabular}{ c @{\hspace{0.25cm}} c @{\hspace{0.25cm}} c @{\hspace{0.25cm}} c}
\includegraphics[scale=0.20]{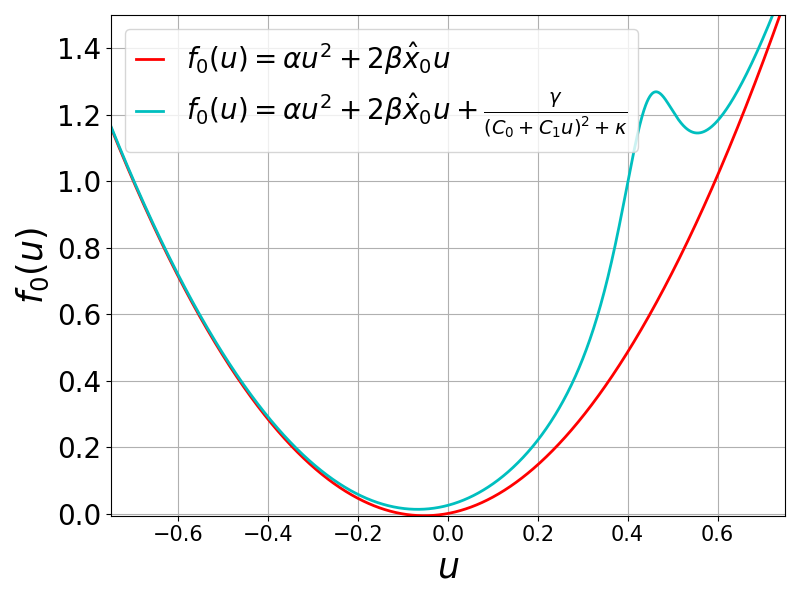} &
\includegraphics[scale=0.20]{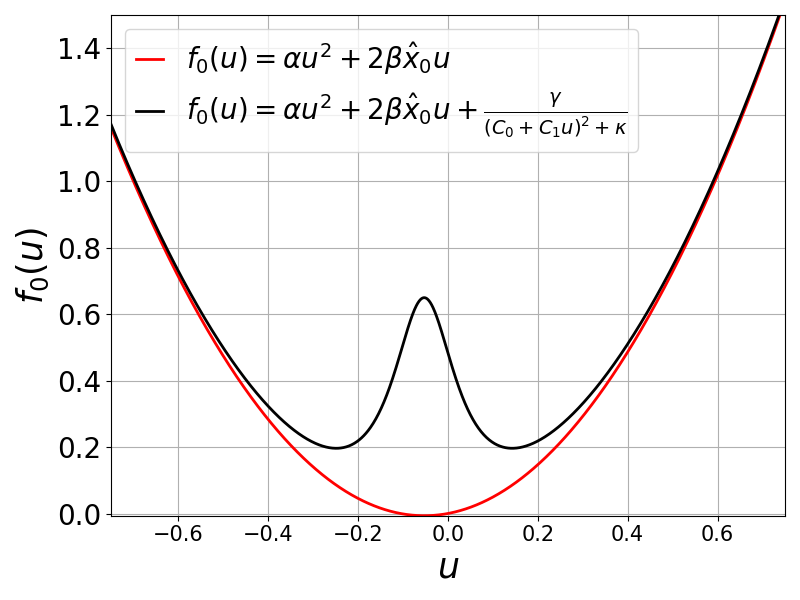} &
\includegraphics[scale=0.20]{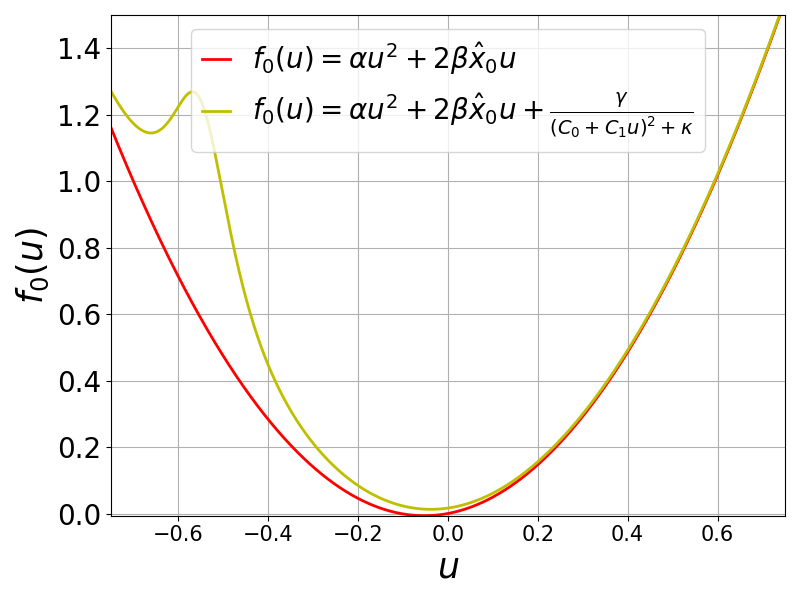} &
\includegraphics[scale=0.20]{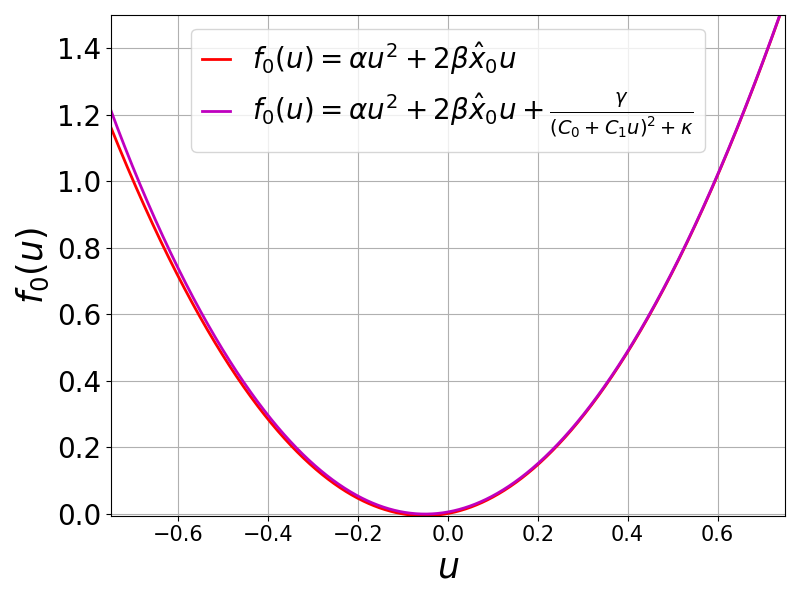} \\
\small \quad (a) $-\frac{C_0}{C_1} = -\frac{\beta \xhat_0}{\alpha}+ 0.5$ &      \small \quad (b) $-\frac{C_0}{C_1} = -\frac{\beta \xhat_0}{\alpha}$ & \small \quad (c) $-\frac{C_0}{C_1} = -\frac{\beta \xhat_0}{\alpha}-0.5$ &      \small \quad(d) $-\frac{C_0}{C_1} = -\frac{\beta \xhat_0}{\alpha}- 1.0$ 
\end{tabular}
\vspace{5pt}
\caption{The landscape of $f_0(u) = f_{LQG}(u) + g(u)$ depends on the distance between the global minima of $f_{LQG}(u):=\alpha u^2 + 2 \beta \xhat_0 u$, and the global maxima of $g(u):=\frac{\gamma}{(C_0+C_1 u)^2 + \kappa}$, where $\alpha, \beta, \gamma, \kappa$ are as defined in \eqref{eqn:alpha beta gamma kappa}. When this distance, captured by $|C_0/C_1 - \beta \xhat_0/\alpha|$ is large $u_0^\star$ and $u_0^{\rm LQG}$ give similar performance. On the other hand, when $|C_0/C_1 - \beta \xhat_0/\alpha|$ is small, $u_0^\star$ and $u_0^{\rm LQG}$ give very different performance.
}
\label{figure1}
\end{center}
\vspace{-12pt}
\end{figure*}

\subsection{Experimental Setups}
We consider one scalar and two vector experimental setups: a system following~\eqref{eqn: dynamics update} with $n=m=p=1$, and $T=2$, double integrator style dynamics, and a system with observation following the ``orthogonal'' condition in Proposition~\ref{prop:perpendicular}. For the vector systems, we consider fixed trajectory length of $T=100$. 

\subsubsection{Scalar Example}
In this experiment, we consider an instance of \eqref{eqn: dynamics update} with $n=m=p=1$, and $T=2$. Specifically, we set $A=0.9$, $B=1$, $Q_T=Q=R=1$. The initial state is sampled from $x_0 {\sim} \Ncal(0.1, 2)$, whereas,  $\{w_0, w_1\}{\distas} \Ncal(0, 0.01)$, and $\{z_0, z_1\}{\distas} \Ncal(0, 0.09)$. 
This gives, $\alpha=2.405$, $\beta =0.126$, $\gamma=0.03$, and $\kappa=0.045$ in \eqref{eqn:alpha beta gamma kappa}. 
$C_1$ is fixed to $C_1=\alpha$, whereas, $C_0$ is varied between $(-3,3)$ such that $|\frac{C_0}{C_1} - \frac{\beta \xhat_0}{\alpha}| \in [0,1]$.
\subsubsection{Double Integrator Style Dynamics}
In this experiment, the dynamics are governed by the following state-space equations:
\begin{align*}
    \vx_{t+1} &= \begin{bmatrix} 1 & h \\ 0 & 1 \end{bmatrix} \vx_t + \begin{bmatrix} 0 \\ h \end{bmatrix}u_t+\vw_t, \\
    y_t &= (C_0+C_1u_t)\begin{bmatrix} 1 & 0 \end{bmatrix} \vx_t + z_t,
\end{align*}
where $\vw_t {\distas} \Ncal(0,0.01\vI_2)$, $z_t {\distas} \Ncal(0,0.01)$, and we use $h=0.3$.
For quadratic cost, we use $\vQ=\vQ_T=\vI_2$, and $\vR=1000$.
The initial state distribution is $\Ncal(0,\vI_2)$.
These dynamics approximate a double integrator with discretization step $h$.
For example, if the state contains position and velocity, then the input is a force.
The bilinear observation model posits that the signal-to-noise ratio of the position sensor scales with the input force.
Put differently, the position measurement results in the application of a force.

\begin{figure*}[t]
\setlength{\abovecaptionskip}{0pt}
\setlength{\belowcaptionskip}{0pt}
\begin{center}
\begin{tabular}{ c @{\hspace{0.25cm}} c @{\hspace{0.25cm}} c @{\hspace{0.25cm}} c}
\includegraphics[scale=0.30]{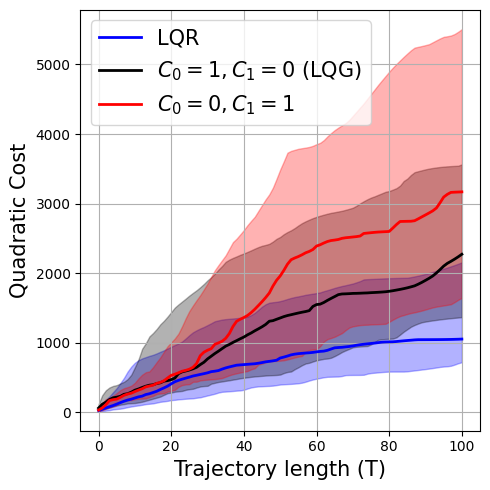} &
\includegraphics[scale=0.30]{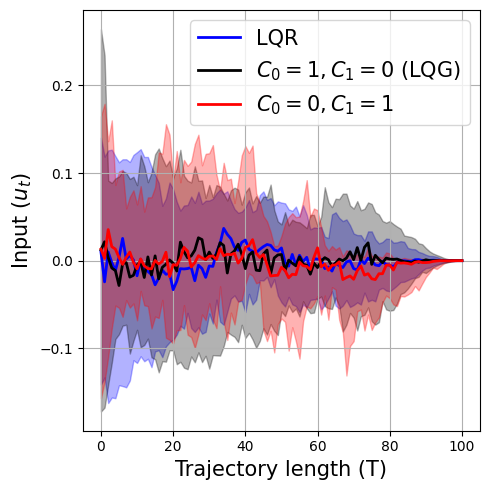} &
\includegraphics[scale=0.30]{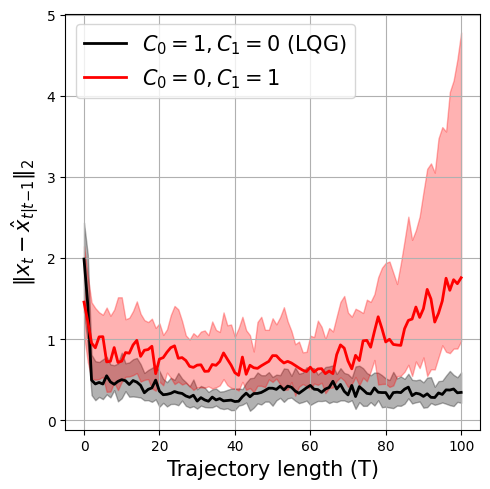} &
\includegraphics[scale=0.30]{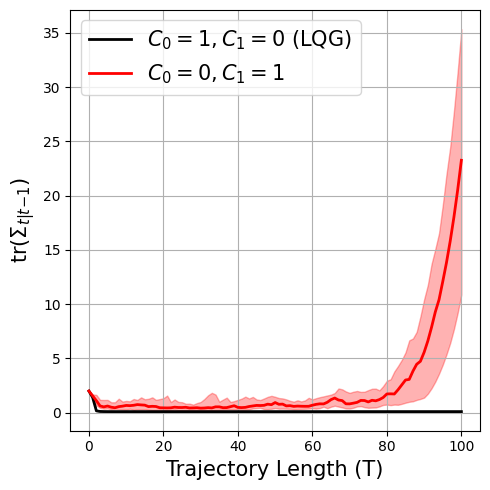} \\
\small \quad (a) Quadratic cost & \small \quad (b) Inputs & \small \quad (c) KF state estimation error & \small \quad(d) KF covariance trace 
\end{tabular}
\vspace{5pt}
\caption{Double Integrator. Bilinear observations $(C_1=1)$ incur more cost (a) than LQG due to the negative effects of small input (b) on state estimation (c,d).}
	\label{figure2}
\end{center}
\vspace{-12pt}
\end{figure*}

\subsubsection{Orthogonal observations}
We consider the system (\ref{eqn: dynamics update}) with $n=6$, and $p=m=3$. The dynamics matrix $\vA$ is generated with i.i.d. $\Ncal(0,1)$ entries, and scaled to have its large eigenvalue $\rho(\vA)=1.1$. The matrix $\vB$ is generated with i.i.d. $\Ncal(0,1/n)$ entries, whereas the matrices $\vC_1,\dots,\vC_p$ are generated with i.i.d. $\Ncal(0,1/m)$ entries. The matrix $\vC_0$ is chosen in the orthogonal complement of the span of $\vC_1,\dots,\vC_p$.
As a result, the system satisfies the condition in Proposition~\ref{prop:perpendicular}.
The process and observation noise follows $\Ncal(0,0.01\vI_n)$ and $\Ncal(0,0.01\vI_m)$, respectively. For quadratic cost, we choose $\vQ=\vQ_T=\vI_n$, and $\vR=\vI_p$. The initial state distribution is $\Ncal(0,\vI_n)$.

\subsection{Scalar Example}
Figure~\ref{figure1} plots $f_0(u)$ defined in \eqref{eqn: f of zero expression} over different values of $u$. Our plots show that the critical points of $f_0(u)$ depends on the distance between the optimal LQG controller $u_0^{\rm LQG} = -\frac{\beta \xhat_0}{\alpha}$ and the global maxima of the estimation error variance term $\frac{\gamma}{(C_0 + C_1 u)^2+\kappa}$, which is $-\frac{C_0}{C_1}$ in this example. Interestingly, Figure~\ref{figure1} shows that when $\frac{C_0}{C_1} = \frac{\beta \xhat_0}{\alpha}$, the LQG controller $u_0^{\rm LQG}=-\frac{\beta \xhat_0}{\alpha}$ locally maximizes $f_0(u)$, which is non-convex, However, when $\frac{C_0}{C_1}$ is far from $\frac{\beta \xhat_0}{\alpha}$, the optimal nonlinear controller converges to the LQG controller.

\subsection{Simulation Results}
For each vector system, we evaluate the performance of the separation principle controller for different observation models.
We consider perfect state observation (LQR), linear state observation (LQG), and bilinear state observation.
For the partially observed models, we use the Kalman filter to estimate the state.
The initial state estimate is sampled at random from the initial distribution to avoid the degenerate case where $\hat{\vx}_t=0$ for all $t\geq0$~(see the discussion preceding Proposition~\ref{prop:perpendicular}).  
We simulate trajectories 50 times and plot the median quantities (solid line) and the 25th and 75th percentiles (shaded).

\subsubsection{Double Integrator}
In Figure~\ref{figure2}(a), we plot the quadratic cost over time, which shows that the bilinear observation model generally incurs higher cost.
Plotting the input over time in Figure~\ref{figure2}(b) reveals that the input magnitude drops towards zero near the end of the trajectory, due to the optimal state feedback policy.
Though this has no ill effect on LQG, it disrupts the state estimation accuracy for the bilinear model, as illustrated in  Figure~\ref{figure2}(c) and (d).
The Kalman filter diverges because the small inputs lead to a loss of observability, highlighting the sub-optimality of the separation principle in this setting.

\begin{figure}[h!]
    \centering
    \begin{subfigure}[h]{0.4\linewidth}
        \centering
        \includegraphics[width=\linewidth]{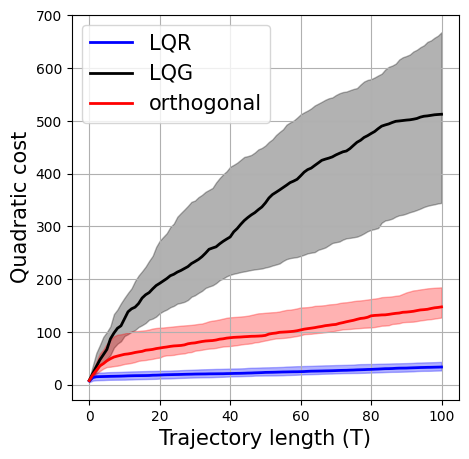}
        \caption{Linear and bilinear observations show different quadratic cost}
\end{subfigure}
    ~
    \begin{subfigure}[h]{0.4\linewidth}
        \centering
        \includegraphics[width=\linewidth]{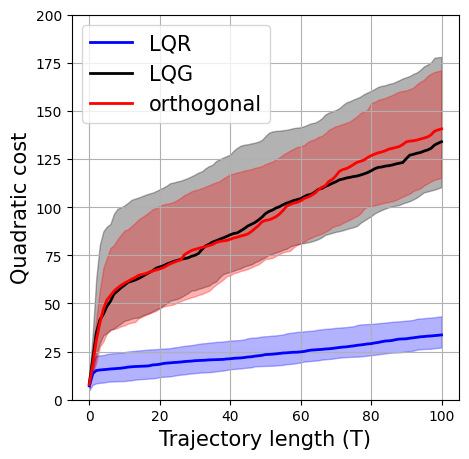}
        \caption{Linear and bilinear observations show similar quadratic cost}
\end{subfigure}
    \caption{For fixed $\vA$, $\vB$, $\vC_1,\dots,\vC_p$, two plots correspond to different $\vC_0$.}
    \label{figure3}
\end{figure}

\subsubsection{Orthogonal Observations}
Following Proposition~\ref{prop:perpendicular}, we consider scenarios where the observation model rules out a loss of observability by design.
For fixed $\vA,\vB, \vC_1,\dots,\vC_p$, we choose two different $\vC_0$ in the orthogonal complement of the span of $\vC_1,\dots,\vC_p$. In Figure~\ref{figure3}, we can see that the choice of $\vC_0$ can show different behavior. Figure~\ref{figure3}(a) shows that the bilinear observations can actually improve the quadratic cost compared with LQG. On the other hand, Figure~\ref{figure3}(b) shows that the performance can be similar for linear vs.\ bilinear observations. 

\section{Conclusion and Discussion}
We study the problem of minimizing a quadratic cost of controlling linear dynamical systems from bilinear observations. 
Our results show that the separation principle does not hold in general for these problems. 
Moreover, the optimal controller is not affine in the estimated state, although the estimation problem itself remains straightforward, that is, the KF still gives the optimal state estimates, and the posterior distribution of the state. 
We find that the optimal control problem has a cost-to-go which is generally nonconvex, and the control inputs affect the estimated state and the estimation error covariance in non-trivial ways.  
We derive analytical expression for the optimal nonlinear control policy in a simple setting and find that the optimal controller is indeed nonlinear in the estimated state. 
We also derive conditions which guarantee uniform observability in the case of bilinear observations, and verify it through numerical experiments.

There are several open questions leading to interesting directions for future work. First, deriving an optimal nonlinear control law for~\eqref{eqn:bilinear LQG} when $T> 2$ is still an open problem. An alternate approach to this will be direct policy optimization for~\eqref{eqn:bilinear LQG}. Second, when the dynamics matrices are unknown, it will be challenging to design an adaptive control scheme for the system~\eqref{eqn: dynamics update}. Lastly, extending this problem to study the optimal control of partially observed bilinear dynamical systems is also an important future direction. 

\section*{Acknowledgements}
S.D. was partly supported by NSF CCF 2312774, NSF OAC-2311521, NSF IIS-2442137, a PCCW Affinito-Stewart Award, a Gift to the LinkedIn-Cornell Bowers CIS Strategic Partnership, and an AI2050 Early Career Fellowship program at Schmidt Sciences. M.F. was supported in part by awards NSF TRIPODS II 2023166, NSF CCF 2007036, NSF CCF 2212261, NSF CCF 2312775, and by an Amazon-UW Hub gift award. Y.J. was partially supported by the Knut and Alice Wallenberg Foundation Postdoctoral Scholarship Program at MIT-KAW 2022.0366.

\newpage 
\bibliographystyle{alpha}
\bibliography{Bibfiles}

\end{document}